\documentclass{amsart}
\usepackage{mathrsfs}
\usepackage{amsmath}
\usepackage{epsfig}
\usepackage{amsfonts}
\usepackage{amssymb}
\usepackage{amsthm}
\usepackage{amscd}
\usepackage[all]{xy}
\usepackage{rotating}
\usepackage{lscape}
\usepackage{amsbsy}
\usepackage{verbatim}
\usepackage{moreverb}

\usepackage{tikz}
\usepackage{calc}

\newcommand{\Vect}{\mathrm{Vect}}
\newcommand{\Hilb}{\mathrm{Hilb}}

\newcommand{\TC}{\mathrm{TC}}
\newcommand{\CN}{\mathrm{CN}}

\newcommand{\Cattens}{\mathrm{Cat}^\otimes}

\newcommand{\Bord}{\mathrm{Bord}}
\newcommand{\Int}{\mathrm{Int}}
\newcommand{\vNalg}{\mathrm{vNalg}}
\newcommand{\OrBord}{\mathrm{OrBord}}
\newcommand{\FrBord}{\mathrm{FrBord}}
\newcommand{\SFBord}{\mathrm{SFBord}}
\newcommand{\ra}{\rightarrow}
\newcommand{\xra}{\xrightarrow}
\newcommand{\cA}{\mathcal{A}}
\newcommand{\cB}{\mathcal{B}}
\newcommand{\cC}{\mathcal{C}}
\newcommand{\cD}{\mathcal{D}}
\newcommand{\cE}{\mathcal{E}}
\newcommand{\cF}{\mathcal{F}}
\newcommand{\cG}{\mathcal{G}}
\newcommand{\cN}{\mathcal{N}}
\newcommand{\cM}{\mathcal{M}}
\newcommand{\cL}{\mathcal{L}}
\newcommand{\cP}{\mathcal{P}}
\newcommand{\fg}{\mathfrak{g}}
\newcommand{\CC}{\mathbb{C}}
\newcommand{\RR}{\mathbb{R}}
\newcommand{\Alg}{\mathrm{Alg}}

\DeclareMathOperator{\Rep}{Rep}
\DeclareMathOperator{\Hom}{Hom}
\DeclareMathOperator{\Dim}{Dim}

\DeclareMathOperator{\Aut}{Aut}
\newcommand{\Nghd}{\mathrm{Nghd}}
\newcommand{\nin}{\notin}
\newcommand{\op}{\text{op}}
\newcommand{\pnt}{\ast}
\newcommand{\fin}{\text{fin}}
\newcommand{\alg}{\text{alg}}
\newcommand{\rev}{\text{rev}}
\newcommand{\Moon}{\mathrm{Moon}}
\newcommand{\Monst}{\mathrm{Monst}}

\newcommand{\bs}{\begin{smallmatrix}}
\newcommand{\es}{\end{smallmatrix}}
\newcommand{\pbs}{\left( \begin{smallmatrix}}
\newcommand{\pes}{\end{smallmatrix} \right)}
\newcommand{\aoa}{\bs \cA \\ \otimes \\ \cA \es}

\newcommand{\cdast}{\circledast}

\newcommand{\tv}{\hat{v}}

\newcommand{\nn}{\nonumber}
\newcommand{\nid}{\noindent}

\newcommand{\cb}{\raisebox{.6ex-.5\height}}

\def\bicolorint{
\bicolorleft
\bicolorright
}

\def\bicolorleft{
\draw[line width=.4pt] (-5,-1) .. controls (-4,2) and (-3,2) .. (-2,1) .. controls (-1,0) .. (0,0);
}
\def\bicolorright{
\draw[line width=1pt] (0,0) .. controls (1,0) .. (2,-.75) .. controls (2.25,-1) and (2.75,-1) .. (3,-.75) .. controls (4.25,.5) and (4.75,.5) .. (5.5,0);
}

\def\bulletsize{.06}
\def\bulletcolor{black!80!gray}
\def\bull(#1,#2){
	\fill [color=\bulletcolor] (#1,#2) circle (\bulletsize);
}

\def\pnt{\hspace{-.8ex}
\raisebox{0.2ex}{
\begin{tikzpicture}[scale=.8]
\bull(0,0); 
\end{tikzpicture}
}\hspace{-.8ex}
}

\def\pospnt{\hspace{-1ex}
\raisebox{-.35ex}{
\begin{tikzpicture}[scale=.8]
\bull(0,0); 
\node at (0.2,0.1) {$\scriptstyle +$};
\end{tikzpicture}
}\hspace{-1.8ex}
}

\def\negpnt{\hspace{-1ex}
\raisebox{-.35ex}{
\begin{tikzpicture}[scale=.8]
\bull(0,0); 
\node at (0.2,0.1) {$\scriptstyle -$};
\end{tikzpicture}
}\hspace{-1.8ex}
}

\def\tworarrow{\hspace{.1cm}{\setlength{\unitlength}{.46mm}\linethickness{.09mm}
\begin{picture}(8,8)(0,0)\qbezier(0,4)(4,7)(8,4)\qbezier(0,1)(4,-2)(8,1)\qbezier(3.5,4)(3.5,3)(3.5,1.5)
\qbezier(4.5,4)(4.5,3)(4.5,1.5)\qbezier(4,0.8)(4.5,1.7)(5.5,2)\qbezier(4,0.8)(3.5,1.7)(2.5,2)
\qbezier(8,1)(7.4,.2)(7.7,-.7)\qbezier(8,1)(7,1)(6.5,1.5)\qbezier(8,4)(7.4,4.8)(7.7,5.7)
                                                \qbezier(8,4)(7,4)(6.5,3.5)
                       \end{picture}\hspace{.1cm}}}

\newtheorem{theorem}{Theorem}
\newtheorem{proposition}[theorem]{Proposition}

\newtheorem{corollary}[theorem]{Corollary}

\theoremstyle{definition}
\newtheorem{definition}{Definition}

\theoremstyle{remark}
\newtheorem{remark}{Remark}
\newtheorem{example}{Example}

\title{Conformal nets and local field theory}

\author{Arthur Bartels}
\address{Westf\"alische Wilhelms-Universit\"at M\"unster, Mathematisches Institut, Einsteinstr. 62, D-48149 M\"unster, DE}
\email{bartelsa@math.uni-muenster.de}
\author{Christopher L. Douglas}
\address{Department of Mathematics, University of California, Berkeley, CA 94720, USA}
\email{cdouglas@math.berkeley.edu}
\author{Andr\'e G. Henriques}
\address{Mathematisch Instituut, Universiteit Utrecht, 3508 TA Utrecht, NL}
\email{a.g.henriques@uu.nl}

\begin{document}

%%%%%%%%%%%%%%%%%%%%%%%%%%%%%%

%%%%%%%%%%%%%%%%%%%%%%%%%%%%%%%%%%%%%%%%%%%%%

\begin{abstract} We describe a coordinate-free notion of conformal nets as a mathematical model of conformal field theory.  We define defects between conformal nets and introduce composition of defects, thereby providing a notion of morphism between conformal field theories.  Altogether we characterize the algebraic structure of the collection of conformal nets as a symmetric monoidal tricategory.  Dualizable objects of this tricategory correspond to conformal-net-valued 3-dimensional local topological quantum field theories.  We prove that the dualizable conformal nets are the finite sums of irreducible nets with finite $\mu$-index.  This classification provides a variety of 3-dimensional local field theories, including local field theories associated to central extensions of the loop groups of the special unitary groups.
\end{abstract}

%\vspace*{-15pt}
\maketitle

Mathematical models of quantum field theory have transformed our understanding of topological invariants of low-dimensional manifolds.  Witten initiated this process in 1989 by introducing a collection of 3-manifold invariants defined as Feynman path integrals of the Chern-Simons action over a space of connections~\cite{witten}.  This formulation of topological invariants gave a physical explanation for the Jones, HOMFLY, and Kaufmann polynomial invariants of knots, and provided new knot and 3-manifold invariants by considering path integrals over spaces of $G$-connections for arbitrary simple Lie groups $G$.  By construction, these 3-manifold invariants came from ``topological quantum field theories" associating a system of algebraic invariants to all 2-dimensional manifolds and 3-dimensional manifolds with boundary.  In 1991, Reshetikhin and Turaev gave a mathematical construction of these systems of invariants, and proved that the invariants come from an ``extended topological quantum field theory" that gives invariants of 1-manifolds, 2-manifolds with boundary, and 3-manifolds with codimension-2 corners~\cite{reshturaev}.  Physical quantum field theories associate algebraic data to neighborhoods of manifolds of any codimension, in particular to regions of spacetime localized to a neighborhood of a single point.  It has been an open question whether there exist such fully ``local" topological quantum field theories, which provide algebraic invariants of 0-manifolds, 1-manifolds with boundary, 2-manifolds with codimension-2 corners, and 3-manifolds with codimension-3 corners.  The purpose of this paper is to outline a proof that there does indeed exist a large class of local topological quantum field theories in dimension 3.  

\section*{Overview and statement of results}
\subsection*{Algebraic structure of conformal field theories}
Our construction of 3-dimensional local field theories relies on a new understanding of the algebraic structure of the collection of all conformal field theories.  Specifically, we introduce a notion of morphism between conformal field theories, and morphisms between these morphisms, and morphisms between these morphisms of morphisms---altogether we recognize that conformal field theories form the objects of a tricategory.  Zero-dimensional manifolds also form the objects of a tricategory, whose 1-, 2-, and 3-morphisms are respectively 1-, 2-, and 3-dimensional manifolds.  Conformal field theories are therefore a natural candidate for the algebraic invariants associated to 0-manifolds in a 3-dimensional local field theory.  

We use a particular mathematical formulation of conformal field theories, namely conformal nets.  The basic data of a conformal net, inspired by algebraic quantum field theory, is the collection of algebras of spacetime local observables.  We define a new notion of coordinate-free, not-necessarily chiral conformal nets, henceforth referred to simply as conformal nets, which permits us to give a mathematical formulation of the notion of defect between conformal nets and, crucially, of the process of fusion of defects.  Altogether, we completely characterize the algebraic structure of the collection of conformal nets: 

\begin{theorem} \label{thm-cn3}
Conformal nets, defects, sectors, and intertwiners form the objects, morphisms, 2-morphisms, and 3-morphisms of a symmetric monoidal tricategory.
\end{theorem} 
\nid Though conceptually straightforward, the precise definition of a symmetric monoidal tricategory is rather involved and will be suppressed in this paper.  Denote the tricategory of conformal nets by $\CN$.  We develop local field theories with values in conformal nets, which is to say symmetric monoidal functors from the tricategory of 0-, 1-, 2-, and 3-manifolds to $\CN$.

\subsection*{Classification of dualizable conformal nets}
Topological quantum field theories have the feature that the invariants associated to large manifolds are determined by gluing together the invariants associated to small manifolds.  In the same vein, a local topological quantum field theory is entirely determined by the algebraic invariant associated to a point.  To produce conformal-net-valued local field theories, we must therefore answer the question, ``When is a conformal net the invariant assigned to a point by a local field theory?"  By a theorem of Hopkins and Lurie~\cite{hopklurie}, the point invariant of an $n$-dimensional local field theory is always a dualizable object of the target algebraic $n$-category, and conversely any dualizable object is the point invariant of a local field theory.  We precisely classify the dualizable conformal nets: 

\begin{theorem}\label{thm-dualizable}
A conformal net is dualizable if and only if it is a finite direct sum of irreducible nets with finite $\mu$-index.
\end{theorem}
\nid Here the $\mu$-index of an irreducible conformal net is the dimension of the representation category of the net, which is to say the sum of the square dimensions of all representations of the net~\cite{klm}.  We refer to nets that are finite sums of irreducible finite $\mu$-index nets as ``finite".  Altogether we have the following: 

\begin{corollary}\label{cor-lft}
For any finite conformal net $\cN$, there exists a 3-dimensional local topological quantum field theory $F$ such that the invariant $F(\pnt)$ assigned to the point $\pnt$ is the net $\cN$.
\end{corollary}
\nid Note that this corollary depends on the compatibility of our preferred notion of symmetric monoidal tricategory with the one utilized by Hopkins and Lurie.

Various conformal nets are known to be finite.  For example, there is a finite net $\cN(\cL G)$ associated to any positive central extension $\cL G$ of the loop group $LG$ of the special unitary group $G=SU(n)$~\cite{wass, xu}---the algebras of observables for this net are the algebras of operators of local loops in the vacuum representation of $\cL G$.  The local field theory $F_{\cN(\cL G)}$ associated to such a net assigns an invariant to the 1-manifold $S^1$, the circle---that invariant $F_{\cN(\cL G)}(S^1)$ is an algebra whose category of modules is expected to be isomorphic as a fusion category to the category of representations of $\cL G$, which is to say to the invariant associated to a circle in the Reshetikhin-Turaev extended field theory.  The nets associated to positive central extensions of other loop groups are also expected to be finite~\cite{gabfrol}, but this has not yet been established.  Another example of a finite net, constructed by Kawahigashi and Longo~\cite{kawalongo}, is the Moonshine net whose automorphism group is the Monster group.

\subsection*{Background and synopsis}

A 3-dimensional TFT, in the sense of Atiyah and Segal, is a symmetric monoidal functor from the bordism category $\Bord_2^3$ of 2-manifolds and 3-manifold bordisms, to the category $\Vect$ of vector spaces.  An extended TFT, a la Reshetikhin and Turaev, is a symmetric monoidal functor from the extended bordism bicategory $\Bord_1^3$ of 1-, 2-, and 3-manifolds, to a target bicategory $\cB$.  Not every bicategory $\cB$ is an appropriate target for an extended TFT: the restriction of the extended TFT functor to the ordinary bordism category $\Bord_2^3$ should produce an ordinary TFT.  The bordism bicategory deloops the ordinary bordism category, in the sense that $\Omega \Bord_1^3 := \Hom_{\Bord_1^3}(1,1) = \Bord_2^3$.  The target bicategory $\cB$ must similarly deloop the category of vector spaces, that is satisfy the condition $\Omega \cB = \Vect$.  Two natural bicategories that have this property are the bicategory of $\Vect$-module categories and the bicategory of algebras; the former is more categorical in character, the latter more algebraic.

A local TFT is a symmetric monoidal trifunctor from the tricategory $\Bord_0^3$ of 0-, 1-, 2-, and 3-manifolds to a target tricategory $\cC$.  The bordism tricategory double deloops the ordinary bordism category, and similarly the target tricategory should satisfy the condition $\Omega^2 \cC = \Vect$.  The tricategory $\TC$ of tensor categories satisfies this condition, and so can function as a target for local field theory.  The conformal net tricategory $\CN$ is a more algebraic, less categorical alternative.  Indeed, we will see that conformal nets have underlying vector spaces with not one, but two distinct algebra structures, and in this sense represent a theory of higher or ``2"-algebra.

In the first section below, we give a precise definition of the notion of 2-algebra.  We then define conformal nets and observe that a conformal net has an underlying 2-algebra.  We introduce defects, fusion of defects, and sectors between defects, and discuss the proof of Theorem~\ref{thm-cn3}.  In the second section, we review the notion of dualizable objects in higher categories, and discuss finiteness properties for nets, defects, and sectors.  We explain the relationship between dualizability and finiteness for nets and describe the proof of Theorem~\ref{thm-dualizable}.  Along the way we construct not just a set but a tricategory of conformal-net-valued local field theories.  We conclude by discussing the loop group nets and their associated field theories.

\section{2-algebras and conformal nets}

\subsection*{2-algebras}
In 1962, Eckmann and Hilton~\cite{eckmannhilton} observed that if $V$ is a vector space with two commuting unital algebra structures, then the two algebra structures are commutative and equal.  Denoting one multiplication horizontally and the other vertically, this follows from the equation
\begin{equation} \nn
a b = \pbs a \\ 1 \pes \pbs 1 \\ b \pes = \bs (a 1) \\ (1 b) \es = \bs a \\ b \es = \bs (1 a) \\ (b 1) \es = \pbs 1 \\ b \pes \pbs a \\ 1 \pes = b a.
\end{equation}
\nid This simple observation has stymied investigation of multiple algebra structures on vector spaces.  We avoid the Eckmann-Hilton dilemma by weakening the associativity condition on the second of the two multiplications, demanding that the vertical multiplication is associative only up to a horizontal conjugation: 
\begin{definition} \label{def-2alg}
A \emph{2-algebra} is a unital associative algebra $\cA$, together with an algebra homomorphism $\mu: \aoa \ra \cA$, denoted $\bs a \\ b \es := \mu \pbs a \\ \otimes \\ b \pes$, and an invertible element $v \in \cA$ such that
\begin{align}
& \hspace*{-185pt} \mathit{1.} \quad v \pbs \pbs a \\ b \pes \\ c \pes v^{-1} = \pbs a \\ \pbs b \\ c \pes \pes , \nn \\
& \hspace*{-185pt} \mathit{2.} \quad v v = \pbs 1 \\ v \pes v \pbs v \\ 1 \pes , \nn
\end{align}
\end{definition}
\nid Notice that because $\mu$ is an algebra homomorphism, the two multiplications of a 2-algebra commute strictly; that is, $\bs (a b) \\ (c d) \es = \pbs a \\ c \pes \pbs b \\ d \pes$, where we have denoted the associative algebra multiplication of $\cA$, as usual, by $a b$.
\begin{remark}
Though our notion of 2-algebra is completely algebraic, the reader may recognize the echos of the homotopical structure of $E_2$-algebras, as follows.  Let $\cA$ be a 2-algebra.  We define $a \in \cA$ and $b \in \cA$ to be \emph{homotopic} if $b = v a v^{-1}$, for $v$ the distinguished element in the definition of a 2-algebra; denote by $\sim$ the equivalence relation multiplicatively generated by this notion of homotopy.  For each element $a \in \cA$, there is an associated element $[a] := \pbs \pbs 1 \\ a \pes \\ 1 \pes$.  Observe that
\begin{equation} \nn
\begin{split}
[a] [b] 
& \sim v [a] v^{-1} [b] 
= v \pbs \pbs 1 \\ a \pes \\ 1 \pes v^{-1} \pbs \pbs 1 \\ b \pes \\ 1 \pes \\
& = \pbs 1 \\ \pbs a \\ 1 \pes \pes \pbs \pbs 1 \\ b \pes \\ 1 \pes 
= \bs \pbs 1 \\ b \pes \\ \pbs a \\ 1 \pes \es 
= \pbs \pbs 1 \\ b \pes \\ 1 \pes \pbs 1 \\ \pbs a \\ 1 \pes \pes \\
& = \pbs \pbs 1 \\ b \pes \\ 1 \pes v \pbs \pbs 1 \\ a \pes \\ 1 \pes v^{-1} 
= [b] v [a] v^{-1} 
\sim [b] [a] .
\end{split}
\end{equation} 
%\enlargethispage{\baselineskip}
In this sense the horizontal multiplication is ``algebraically homotopically commutative".
\end{remark} 
\begin{remark}
In a monoidal category $\cC$, the tensor map $\otimes: \cC \times \cC \ra \cC$ is not strictly associative, but is subject to the existence of an associator transformation $\cC \times \cC \times \cC \!\!\raisebox{-.15ex}{\tworarrow} \cC$ satisfying the pentagon axiom.  The distinguished element $v \in \cA$ in a 2-algebra plays a role analogous to the associator transformation: by condition (1) of Definition~\ref{def-2alg}, horizontal conjugation by $v$ implements vertical associativity.  Condition (2) on $v$, namely $v v = \pbs 1 \\ v \pes v \pbs v \\ 1 \pes$, is therefore a purely algebraic analogue of the categorical pentagon axiom.
\end{remark}
%\enlargethispage{\baselineskip}

\subsection*{Conformal nets}
Conformal nets are, in a sense we make precise presently, topological 2-algebras.  Conformal nets are usually defined as nets of von Neumann algebras on the standard circle.  Recall that von Neumann algebras are unital $\ast$-subalgebras of the bounded operators on complex Hilbert space that are closed in the weak topology, or alternately in the topology of pointwise convergence.  We adopt a more coordinate-free perspective on nets, as follows.  Let $\Int$ denote the category whose objects are oriented intervals, that is non-empty, connected, simply connected, compact 1-manifolds with boundary, and whose morphisms are not-necessarily orientation-preserving embeddings.  Let $\vNalg$ denote the category whose objects are von Neumann algebras and whose morphisms are all homomorphisms and antihomomorphisms.
\begin{definition}
A \emph{conformal net} $\cN$ is a cosheaf of von Neumann algebras on the category of intervals.  That is, a conformal net is a continuous covariant functor $\cN\!:\! \Int \!\ra\! \vNalg$ taking orientation-preserving embeddings to homomorphisms and orientation-reversing embeddings to antihomomorphisms, subject to the following conditions, for $\!I\!$ and $\!J\!$ subintervals of the interval $\!K\!$:
\begin{enumerate}
\item[1.] \textbf{Isotony} The induced map $\cN(I) \ra \cN(K)$ is injective. \vspace{2pt}
\item[2.] \textbf{Additivity} $\cN(I \cup J) = \cN(I) \vee \cN(J)$; here $\vee$ indicates the von Neumann algebra generated by the two algebras. \vspace{2pt}
\item[3.] \textbf{Locality} If $I \cap J$ is a point, then $\cN(I)$ and $\cN(J)$ are commuting subalgebras of $\cN(K)$.\vspace{2pt}
\item[4.] \textbf{Split property} If $I \cap J = \emptyset$ then the map $\cN(I) \otimes^{\text{alg}} \cN(J) \ra \cN(K)$ extends to the spatial tensor product $\cN(I) \bar{\otimes} \cN(J)$. \vspace{2pt}
\item[5.] \textbf{Inner covariance} If $\phi: I \ra I$ is a diffeomorphism that is the identity in a neighborhood of the boundary, then $\cN(\phi)$ is an inner automorphism. \vspace{2pt}
\item[6.] \textbf{Vacuum sector} Suppose $I \subsetneq K$ contains the boundary point $p \in \partial K$, and $\bar{I}$ denotes $I$ with the reversed orientation; the algebra $\cN(I)$ acts on $L^2(\cN(K))$ via the left action of $\cN(K)$, and $\cN(\bar{I})$ acts on $L^2(\cN(K))$ via the right action of $\cN(K)$.  In this case, the action of $\cN(I) \otimes^{\text{alg}} \cN(\bar{I})$ on $L^2(\cN(K))$ extends to an action of $\cN(I \cup_p \bar{I})$. \vspace{2pt}
\end{enumerate}
\end{definition}
\nid A net is irreducible, by definition, if for all $I$ the algebra $\cN(I)$ is a factor, that is the center $Z(\cN(I))$ is $\CC$.  Note that, as all intervals are isomorphic, a net $\cN$ has a uniquely determined underlying algebra, namely $\cN([0,1])$.  The definition of a net is partially elucidated by observing that this algebra is in fact a 2-algebra:
\begin{proposition}
For $\cN$ a conformal net, the algebra underlying the von Neumann algebra $\cN([0,1])$ is a 2-algebra.
\end{proposition}
\nid The second algebra structure on $\cN([0,1])$ comes from concatenation of intervals; the interval $[0,1]$ should be pictured vertically, to indicate that interval concatenation gives the weakly associative vertical multiplication in the 2-algebra.  We give the proof under the assumption that the net is irreducible.

\vspace{3pt}
\noindent \emph{Proof:} %\begin{proof}
Let $i: [0,1] \ra [0,2]$ and $j: [0,1] \ra [0,2]$ be respectively the inclusion and the inclusion by unit translation, and let $h: [0,1] \ra [1,2]$ and $k: [1,2] \ra [2,3]$ be the translations and $l: [0,1] \ra [0,1]$ the identity.  Let $s: [0,2] \ra [0,1]$ be a diffeomorphism that has derivative $1$ in a neighborhood of the boundary and let $t: [1,3] \ra [1,2]$ denote the translation of $s$.  The homomorphism $\mu: \bs \cN([0,1]) \\ \otimes \\ \cN([0,1]) \es \ra \cN([0,1])$ is defined by
\begin{equation} \nn
\mu \pbs a \\ \otimes \\ b \pes := s_* (j_*(a) i_*(b)).
\end{equation}

Let $p: [0,1] \ra [0,1]$ be the diffeomorphism $s (s \cup k^{-1}) (l \cup t^{-1}) s^{-1}$, and let $\pbs 1 \\ p \pes: [0,1] \ra [0,1]$ and $\pbs p \\ 1 \pes: [0,1] \ra [0,1]$ denote respectively the diffeomorphisms $s (p \cup (h l h^{-1})) s^{-1}$ and $s (l \cup (h p h^{-1})) s^{-1}$.  For any choice of diffeomorphism $s$ we have $p p = \pbs 1 \\ p \pes p \pbs p \\ 1 \pes$.  Let $\tv \in \cN([0,1])$ be an element implementing the inner automorphism $p_*$; that is $\tv a \tv^{-1} = p_* a$.  Such an element exists by the covariance property of $\cN$.  The additivity property of $\cN$ in turn implies that the elements $\pbs 1 \\ \tv \pes$ and $\pbs \tv \\ 1 \pes$ implement respectively the automorphisms $\pbs 1 \\ p \pes_*$ and $\pbs p \\ 1 \pes_*$.  The product $\lambda := \pbs \tv^{-1} \\ 1 \pes \tv^{-1} \pbs 1 \\ \tv^{-1} \pes \tv \tv$ is in the center of $\cN([0,1])$ and is therefore scalar.  Define the distinguished element $v \in \cN([0,1])$ by $v := \lambda \tv$ and note that $vv = \pbs 1 \\ v \pes v \pbs v \\ 1 \pes$, as desired.  Moreover, by construction the element $v$ satisfies the condition
$v \pbs \pbs a \\ b \pes \\ c \pes v^{-1} = p_* \pbs \pbs a \\ b \pes \\ c \pes = \pbs a \\ \pbs b \\ c \pes \pes$.
%\end{proof}

\begin{remark}
For technical reasons involved in the proof of Theorem~\ref{thm-111} below, we restrict attention to conformal nets that are direct sums of irreducible nets with finite $\mu$-index---see the next section for a discussion of the $\mu$-index.  Henceforth, we let `conformal net' refer to only nets of this form.
\end{remark}

Algebras are the objects of a symmetric monoidal bicategory; we might expect 2-algebras to form the objects of a symmetric monoidal tricategory.  We prove that this is indeed the case, provided we work not with 2-algebras but with conformal nets.  
The loop category of algebras is the category of vector spaces, and the double loop category of conformal nets should also be a category of vector spaces---in fact the objects of this double loop category will have the structure of Hilbert spaces.  Let $\Cattens$ denote the 2-category of symmetric monoidal categories. The symmetric monoidal bicategory of algebras is conveniently encoded as a category internal to $\Cattens$: the object symmetric monoidal category is algebras and algebra homomorphisms, and the morphism symmetric monoidal category is bimodules and bimodule maps.  For conformal nets, we analogously encode the structure of a symmetric monoidal tricategory as a bicategory internal to $\Cattens$.  The following is therefore a slight refinement of Theorem~\ref{thm-cn3}:
\begin{theorem} \label{thm-cnbicat}
Conformal nets, defects, and sectors between defects form a bicategory $\CN$ in symmetric monoidal categories, such that the double loop category $\Omega^2 \CN$ is the symmetric monoidal category $\Hilb$ of Hilbert spaces.
\end{theorem}
\nid Before we can describe the proof of Theorem~\ref{thm-cnbicat}, we need, of course, to define defects and sectors.  

\subsection*{Defects and sectors}
A conformal field theory associates, in particular, evolution operators to conformal surfaces.  Traditionally, a defect associates operators to a conformal surface which is `bicolored', that is in which some regions are labeled with one CFT and other regions are labeled by a second CFT---the defect resides along the `color' change between these regions and is viewed as a transformation between the two CFTs.  In the formalism of conformal nets, a conformal field theory is described by algebras of operators associated to intervals, and correspondingly a defect is encoded in a collection of algebras of operators associated to bicolored intervals.  By definition a \emph{bicolored} interval is an oriented interval $I$ equipped with two subintervals $I_w$ (the white subinterval) and $I_b$ (the black subinterval) such that $I = I_w \cup I_b$ and either $I_b = \emptyset$ or $I_w = \emptyset$ or $I = I_w \vee I_b$; here $\vee$ indicates the union along a point.  In the last case, the bicolored interval is moreover equipped with a local coordinate $c: \Nghd(I_w \cap I_b) \ra \RR$ around $I_w \cap I_b$, with $c^{-1}(\RR_{\leq 0}) \subset I_w$ and $c^{-1}(\RR_{\geq 0}) \subset I_b$.  We refer to the three kinds of bicolored intervals as white, black, and genuinely bicolored, respectively.  The genuinely bicolored intervals are denoted pictorially as
\begin{tikzpicture}[scale=.1]
\bicolorint
\end{tikzpicture}
.  A morphism of bicolored intervals $I \ra J$ is an embedding $\phi: I \ra J$ of the underlying intervals that preserves the local coordinate and such that $\phi^{-1}(J_w) = I_w$ and $\phi^{-1}(J_b) = I_b$.
\begin{definition}
A \emph{defect} $\cD: \cN \ra \cM$ from the net $\cN$ to the net $\cM$ is a continuous covariant functor $\cD: \Int_{bicol} \ra \vNalg$ from bicolored intervals to von Neumann algebras, such that the restriction of $\cD$ to the full subcategory of white, respectively black, intervals in $\Int_{bicol}$ is naturally isomorphic to the net $\cN$, respectively $\cM$.  The functor should satisfy the following conditions, for $I$ a genuinely bicolored subinterval of and $J$ a bicolored subinterval of the genuinely bicolored interval $K$:
\begin{enumerate}
\item[1.] \textbf{Isotony} The map $\cD(I) \ra \cD(K)$ is injective. \vspace{2pt}
\item[2.] \textbf{Additivity} $\cD(I \cup J) = \cD(I) \vee \cD(J)$. \vspace{2pt}
\item[3.] \textbf{Haag duality} If $I\cap J$ is a point and $K = I \cup J$, then $\cD(I) = \cD(J)'$; here $\cD(J)'$ denotes the commutant of $\cD(J)$ in $\cD(K)$. \vspace{2pt}
\item[4.] \textbf{Vacuum sector} Suppose $J$ is either a white or black interval containing exactly one of the boundary points, say $p$, of the genuinely bicolored interval $K$; let $\bar{J}$ denote $J$ with the reversed orientation.  The action of $\cD(J) \otimes^{\text{alg}} \cD(\bar{J})$ on $L^2(\cD(K))$, via the left and right actions of $\cD(K)$ on $L^2(\cD(K))$, extends to an action of $\cD(J \cup_p \bar{J})$. \vspace{2pt}
\end{enumerate}
\end{definition}

\noindent A defect is irreducible, by definition, if $\cD(I)$ is a factor for all bicolored intervals $I$.
Pictorially, we abbreviate a defect $\cD: \cN \ra \cM$ as 
\begin{tikzpicture}[scale=.1]
\bicolorint
\node[above] at (0,0) {$\scriptstyle \cD$};
\end{tikzpicture}
, where the thinner line indicates the `color' of the net $\cN$ and the thicker line indicates the `color' of the net $\cM$.  We also denote the defect as $_\cN \cD_\cM$ to indicate that the defect functions as a type of higher bimodule between the two nets.  
\begin{example}
Let $\cN\hookrightarrow \cM$ be a non-local extension of a conformal net $\cN$~\cite{longorehren}.
The functor sending $I$ to $\cM(I)$ if $I$ is genuinely bicolored and to $\cN(I)$ otherwise is a defect from $\cN$ to itself.
\end{example}
\begin{example}
Let $\cN\hookrightarrow \cM$ be an inclusion of loop group nets coming from a conformal inclusion---see Bais and Bouwknegt~\cite{bb} for a classification.
The functor sending $I$ to $\cN(I)$ if $I$ is white and to $\cM(I)$ otherwise is a defect from $\cN$ to $\cM$.
\end{example}
It has heretofore been a mystery how to compose arbitrary defects between conformal field theories.  Conceiving of defects as higher bimodules and exploiting the technical framework of coordinate-free conformal nets, we introduce the crucial new operation of composition of defects:
\begin{equation} \nn
(_\cN \cD_\cM) \circ (_\cM \cE_\cL) := \cD \, \cdast_{\!\cM} \, \cE.
\end{equation}
Here $\cdast_{\!\cM}$ is a fusion product over the conformal net $\cM$---this novel construction is described in the proof of Theorem~\ref{thm-cnbicat} below.

Next we define sectors between defects.  Let $S^1$ denote the standard unit circle in $\CC$.  Let $I \subset S^1$ be an interval with $i \nin \partial I$ and $-i \nin \partial I$, and either $i \nin I$ or $-i \nin I$.  Equip any such $I$ with a bicoloring by letting $I_w$ be the subinterval of points with nonpositive real part, and $I_b$ be the subinterval of points with nonnegative real part.
\begin{definition}
Given two defects $_\cN \cD_\cM$ and $_\cN \cE_\cM$ from the net $\cN$ to the net $\cM$, a \emph{sector} $H: \cD \ra \cE$ from $\cD$ to $\cE$ is a Hilbert space $H$, often denoted $_\cD H_\cE$, together with homomorphisms
\begin{eqnarray}
\rho_I : \cN(I) \!\!\!\!&\ra&\!\!\!\! B(H), \text{ for } I \text{ white,} \nn \\
\rho_I : \cM(I) \!\!\!\!&\ra&\!\!\!\! B(H), \text{ for } I \text{ black,} \nn \\
\rho_I : \cD(I) \!\!\!\!&\ra&\!\!\!\! B(H), \text{ for } i \in I, \nn \\
\rho_I : \cE(I) \!\!\!\!&\ra&\!\!\!\! B(H), \text{ for } -i \in I, \nn
\end{eqnarray}
for each interval $I \subset S^1$ with $i \nin \partial I$, $-i \nin \partial I$, and either $i \nin I$ or $-i \nin I$.  These homomorphisms should be compatible in the sense that $\rho_K |_I = \rho_I$ for $I \subset K$, and local in the sense that the images of $\rho_I$ and $\rho_J$ should commute if $I$ and $J$ have disjoint interiors.
\end{definition}
\nid Pictorially, we depict a sector with its actions of the defects and nets as 
\cb{
\begin{tikzpicture}[scale=.3]
\draw[line width=.4pt] (0,0) .. controls (-.555,0) and (-1,.445) .. (-1,1) .. controls (-1,1.555) and (-.555,2) .. (0,2);
\draw[line width=1pt] (0,0) .. controls (.555,0) and (1,.445) .. (1,1) .. controls (1,1.555) and (.555,2) .. (0,2);
\node at (0,1) {$\scriptstyle H$};
\node[left] at (-1,1) {$\scriptstyle \cN$};
\node[right] at (1,1) {$\scriptstyle \cM$};
\node[above] at (0,2) {$\scriptstyle \cD$};
\node[below] at (0,0) {$\scriptstyle \cE$};
\end{tikzpicture}
}
.  A map of sectors is a map of Hilbert spaces $H \ra K$ that is equivariant with respect all the actions $\rho_I$ in the above definition.

\begin{remark}
In some bicategories, there are two natural notions of equivalence, namely isomorphism and Morita equivalence---for the latter, objects $A$ and $B$ are equivalent if there are morphisms $\phi: A \ra B$ and $\psi: B \ra A$ such that $\phi \psi$ and $\psi \phi$ are isomorphic to the identity.  For instance, Morita equivalence of algebras is a crucial tool in representation theory.  There is an analogous notion of `higher Morita equivalence' in a tricategory, namely two objects $N$ and $M$ are equivalent if there are morphisms $D: N \ra M$ and $E: M \ra N$ and morphisms $\phi: E D \ra 1_N$, $\bar{\phi}: 1_N \ra E D$, $\psi: D E \ra 1_M$, and $\bar{\psi}: 1_M \ra D E$ such that $\phi \bar{\phi}$, $\bar{\phi} \phi$, $\psi \bar{\psi}$, and $\bar{\psi} \psi$ are all isomorphic to the identity.
%}
%
%\emph{
Conformal nets are a model for conformal field theories, and so Theorem~\ref{thm-cnbicat}, which establishes that conformal nets form a tricategory, provides a new notion of equivalence, namely higher Morita equivalence, between CFTs.  This notion seems to be unknown even in the physics literature.  
\end{remark}

\subsection*{Proof of Theorem~\ref{thm-cnbicat}}
We describe some of the most important data of the conformal net bicategory $\CN$ internal to (symmetric monoidal) categories.  We first need a category $\CN_0$ of objects, a category $\CN_1$ of morphisms, and a category $\CN_2$ of 2-morphisms. Conformal nets are functors $\Int \ra \vNalg$, and $\CN_0$ is the category of conformal nets and their natural transformations.  Similarly, defects are functors $\Int_{bicol} \ra \vNalg$, and $\CN_1$ is the category of defects and their natural transformations.  Finally, $\CN_2$ is the category of sectors and maps of sectors.

The identity functor $1: \CN_0 \ra \CN_1$ is given on objects as follows.  The identity defect $1_\cN$ from $\cN$ to $\cN$ evaluates on any bicolored interval $I$ as $1_\cN(I) := \cN(I)$---that is, simply forget the bicoloring.  The vertical identity functor $1: \CN_1 \ra \CN_2$ is given on objects as follows.  Let $S^1_+$ and $S^1_-$ denote the upper and lower halves of the standard circle.  The Hilbert space of the identity sector $1_\cD : \cD \ra \cD$ is $L^2(\cD(S^1_+))$, which we sometimes abbreviate as $L^2(\cD)$.  The actions of $\cD(I)$ for $I$ in $S^1_+$ or $S^1_-$ are inherited from the left and right actions, respectively, of $\cD(S^1_+)$ on $L^2(\cD(S^1_+))$.

A vital piece of data in the internal bicategory $\CN$ is the composition of 1-morphisms, which is to say the composition of defects; this is a functor $\cdast: \CN_1 \times_{\CN_0} \CN_1 \ra \CN_1$ given on objects as follows.  For any genuinely bicolored interval $I$, define an associated bicolored interval $[I_w] := I_w \cup [0,1]$ with $[I_w]_w = I_w$ and an associated bicolored interval $[I_b] := I_b \cup [-1,0]$ with $[I_b]_b = I_b$.  These intervals are versions of $I$ with respectively standard black and standard white segments.  Let $J:=[-1,-1/2]$ be considered as a subinterval of $[I_b]$ via inclusion and as a subinterval of $[I_w]$ via $x \mapsto -x$.  Given nets $\cM$, $\cN$, and $\cP$ and defects $_\cM \cD_\cN$ and $_\cN \cE_\cP$, the composition of defects, denoted $\cD \cdast_\cN \cE$, is defined on a genuinely bicolored interval $I$ by
\begin{equation} \nn
(\cD \cdast_\cN \cE) (I) := \cD([I_w]) \cdast_{\cN(J)} \cE([I_b]).
\end{equation}
Pictorially this equation becomes 
\begin{equation} \nn
(\cD \cdast_\cN \cE) (
\begin{tikzpicture}[scale=.1]
\bicolorint
\end{tikzpicture}
) = \cD(
\raisebox{-.75ex}{
\begin{tikzpicture}[scale=.1]
\bicolorleft
\draw[line width=.7pt, color=white!60!black] (0,0) -- (1.8,0) .. controls (2.3,0) and (2.5,-.2) .. (2.5,-.7) -- (2.5,-2.5);
\end{tikzpicture}
}
) \cdast_{\cN(
\raisebox{-.4ex}{
\begin{tikzpicture}[scale=.1]
\draw[line width=.7pt, color=white!60!black] (0,-.4) -- (0,-2.5);
\end{tikzpicture}
}
)} \cE(
\raisebox{-.75ex}{
\begin{tikzpicture}[scale=.1]
\bicolorright
\draw[line width=.7pt, color=white!60!black] (0,0) -- (-1.8,0) .. controls (-2.3,0) and (-2.5,-.2) .. (-2.5,-.7) -- (-2.5,-2.5);
\end{tikzpicture}
}
).
\end{equation}
Here for $U$, $V$, and $W$ von Neumann algebras, with $V^{\op} \ra U$ and $V \ra W$, the expression $U \cdast_V W$ denotes the fiber product von Neumann algebra.  By definition
\begin{equation} \nn
U \cdast_V W := (U'^{B(H)} \otimes W'^{B(K)})'^{B(H \boxtimes_V K)}
\end{equation}
where $H$ and $K$ are faithful left representations of $U$ and $W$ respectively; here $\boxtimes$ denotes Connes fusion and $A'^B$ denotes the commutant of $A$ in $B$.

Next we describe the horizontal and vertical composition of 2-morphisms.  Given nets $\cM$, $\cN$, and $\cP$, defects $_\cM \cD_\cN$, $_\cM \cE_\cN$, $_\cN \cF_\cP$, and $_\cN \cG_\cP$, and sectors $_\cD H_\cE$ and $_\cF K_\cG$, the horizontal composition of $H$ and $K$, denoted $H \circ_h K$ or $H \boxtimes_\cN K$, is a sector from $\cD \cdast_\cN \cF$ to $\cE \cdast_\cN \cG$ defined as
\begin{equation} \nn
H \circ_h K := H \boxtimes_{\cN(I)} K.
\end{equation}
Here $I$ is the interval $[-\pi / 3, \pi / 3]$; the algebra $\cN(I)$ acts on $H$ via the map $I \ra S^1, x \mapsto \exp x$ and acts on $K$ via the map $I \ra S^1, x \mapsto -\exp(-x)$.  Next, given nets $\cM$ and $\cN$, defects $_\cM \cD_\cN$, $_\cM \cE_\cN$, and $_\cM \cF_\cN$, and sectors $_\cD H_\cE$ and $_\cE K_\cF$, the vertical composition of $H$ and $K$, denoted $H \circ_v K$ or  $H \boxtimes_\cE K$, is a sector from $\cD$ to $\cF$ defined as
\begin{equation} \nn
H \circ_v K := H \boxtimes_{\cE(J)} K.
\end{equation}
Here $J$ is the interval $[\pi,2\pi]$; the algebra $\cE(J)$ acts on $H$ via the map $J \ra S^1, x \mapsto \exp x$ and on $K$ via the map $J \ra S^1, x \mapsto \exp (-x)$. 

The most important piece of structure in a tricategory is the commutator between horizontal and vertical composition of 2-morphisms.  In particular, for sectors $_\cD H_{\cD'}$, $_{\cD'} L_{\cD''}$, $_\cE K_{\cE'}$, and $_{\cE'} M_{\cE''}$, we need an isomorphism between $(H \circ_h K) \circ_v (L \circ_h M)$ and $(H \circ_v L) \circ_h (K \circ_v M)$.  The construction of such an isomorphism can be reduced to the case when all the sectors $H$, $L$, $K$, and $M$ are identity sectors, that is $H = L = 1_\cD := L^2(\cD)$ and $K = L = 1_\cE := L^2(\cE)$.  The construction of the isomorphism for identity sectors is a consequence of the following crucial result:
\begin{theorem} \label{thm-111}
Given conformal nets $\cM$, $\cN$, and $\cP$, and defects $_\cM \cD_\cN$ and $_\cN \cE_\cP$, there exists a natural isomorphism of sectors
\begin{equation} \nn
L^2(\cD \cdast_\cN \cE) \cong L^2(\cD) \boxtimes_{\cN} L^2(\cE).
\end{equation}
\end{theorem}

% [end proof]

\section{Dualizable nets and local field theories}

\subsection*{Classification of local field theories}
A 1-dimensional oriented topological field theory, that is a symmetric monoidal functor $\OrBord_0^1 \ra \Vect$, associates a vector space $V_+$ to the positively oriented point
\!\pospnt\!, and a vector space $V_-$ to the negatively oriented point \negpnt.  These vector spaces are dual because there are bordisms
\cb{
\begin{tikzpicture}[scale=.7]
\draw (0,0) .. controls (1,0) and (1,-.5) .. (0,-.5);
\node at (-.15,.1) {$\scriptstyle +$};
\node at (-.15,-0.4) {$\scriptstyle -$};
\end{tikzpicture}
}
and
\cb{
\begin{tikzpicture}[scale=.7, xscale=-1]
\draw (0,0) .. controls (1,0) and (1,-.5) .. (0,-.5);
\node at (-.15,.1) {$\scriptstyle -$};
\node at (-.15,-0.4) {$\scriptstyle +$};
\end{tikzpicture}
}, thus maps $V_+ \otimes V_- \ra \CC$ and $\CC \ra V_- \otimes V_+$ satisfying the relations $
\cb{
\begin{tikzpicture}[scale=.7,yscale=.5]
\draw (0,0) .. controls (1.25,0) and (1.25,-.25) .. (.75,-.5) .. controls (.25,-.75) and (.25,-1) .. (1.5,-1);
\end{tikzpicture}
}
= 
\cb{
\tikz \draw[scale=.7] (0,0) -- (1,0);
}
$ and 
$
\cb{
\begin{tikzpicture}[scale=.7,yscale=-.5]
\draw (0,0) .. controls (1.25,0) and (1.25,-.25) .. (.75,-.5) .. controls (.25,-.75) and (.25,-1) .. (1.5,-1);
\end{tikzpicture}
}
= 
\cb{
\tikz \draw[scale=.7] (0,0) -- (1,0);
}
$.  Conversely, any dualizable vector space $V_+$ is the image of the positively oriented point for a 1-dimensional oriented TFT.  Hopkins and Lurie recognized that, replacing the orientation structure by a framed structure, the same is true in all dimensions:
\begin{theorem}\emph{(\cite{hopklurie})} \label{thm-hl}
Given an $n$-dimensional framed local field theory with target $\cC$, that is a symmetric monoidal functor $F: \FrBord_0^n \ra \cC$, the image of the positively framed point is a dualizable object $F(\,\pospnt\,) \in \cC$, and conversely any dualizable object of $\cC$ is the image of the positively framed point under some framed local field theory.
\end{theorem}
\begin{remark}
Hopkins and Lurie utilize a notion of symmetric monoidal $n$-category based on $n$-fold complete Segal spaces.  Henceforth we operate under the assumption that when $n=3$ that notion is comparable to our preferred notion of symmetric monoidal tricategory, namely bicategory in symmetric monoidal categories.
\end{remark}

We recall the definition of a dualizable object in a symmetric monoidal $n$-category.  An $n$-category $\cC$ is said to have adjoints for 1-morphisms if all 1-morphisms admit left and right adjoints in the homotopy 2-category of $\cC$.  An $n$-category $\cC$ has adjoints for $k$-morphisms, $1 < k < n$, if all Hom $(n-1)$-categories $\Hom_\cC(a,b)$ have adjoints for $(k-1)$-morphisms.  Furthermore, a symmetric monoidal $n$-category $\cC$ \emph{has duals} if (1) every object has a dual in the homotopy category and (2) $\cC$ has adjoints for all $k$-morphisms for $0<k<n$.  An object of a symmetric monoidal $n$-category $\cC$ is \emph{dualizable} if it is in the maximal symmetric monoidal sub-$n$-category of $\cC$ that has duals.  

In light of the above relationship between local field theories and dualizability, it is a fundamental problem to classify the dualizable objects of target categories for field theory.  This classification is known in dimensions 1 and 2: the dualizable objects of $\Vect$ are the finite-dimensional vector spaces, and the dualizable objects of $\Alg = B\Vect$ are the algebras $A$ that are finite-dimensional as vector spaces and projective as $(A \otimes A^{\op})$-modules.  In this section we give the analogous classification in dimension 3 by classifying the dualizable conformal nets, that is objects of $\CN = B^2 \Vect$.  

\subsection*{Finite nets are dualizable} 
As in the case of vector spaces and algebras, the dualizable nets are characterized by a strong finiteness property.  For nets, the relevant invariant is a representation-theoretic dimension called $\mu$-index.  Let $J_1=\exp([-\frac{\pi i}{4},\frac{\pi i}{4}])$, $J_2=\exp([\frac{\pi i}{4},\frac{3\pi i}{4}])$, $J_3=\exp([\frac{3\pi i}{4},\frac{5\pi i}{4}])$, and $J_4=\exp([\frac{5\pi i}{4},\frac{7\pi i}{4}])$; as before, $S^1_+ = \exp([0,\pi])$ and $S^1_- = \exp([\pi,2 \pi])$.
\begin{definition}
The \emph{$\mu$-index} $\mu(\cN)$ of an irreducible conformal net $\cN$ is the square of the statistical dimension of the bimodule
\begin{equation} \nn
_{\cN(J_1) \bar{\otimes} \cN(J_3)} L^2(\cN(S^1_+))_{\cN(J_2)^{\op} \bar{\otimes} \cN(J_4)^{\op}} .
\end{equation}
\end{definition}
\nid Kawahigashi, Longo, and M\"uger proved that
%\begin{equation} \nn
$\mu(\cN) = \Dim(\Rep(\cN))$
%\end{equation}
for irreducible conformal nets in their (somewhat different than our) sense~\cite{klm}.  Here $\Dim(\Rep(\cN))$ is by definition the sum of the squared categorical dimensions of the irreducible representations of $\cN$.  We expect the same equality is true in our context, and thus may think of an irreducible net with finite $\mu$-index as, in the appropriate sense, finite dimensional.  In general we call a conformal net \emph{finite} if it is a finite direct sum of irreducible nets with finite $\mu$-index.

We now describe a collection of finite conformal nets.  Implicitly we use the fact that a net can be defined by specifying the cosheaf of algebras on the skeletal subcategory $\Int_{S^1} \subset \Int$ of subintervals of $S^1$.

\begin{example}[Loop group nets \cite{gabfrol, wass}] \label{eg-lg}
For $G=SU(n)$, let $\pi: \cL G \ra LG$ be the level $k$ central extension of the loop group $LG$.  For $I \subset S^1$, let $\cL_I G$ be the preimage under $\pi$ of the loops supported in the interval $I$.  Let $\rho: \cL G \ra B(H_0)$ be the vacuum representation, and define the net $\cN_{G,k}$ as 
\begin{equation} \nn
\cN_{G,k}(I) := \langle \rho(\ell) \vert \ell \in \cL_I G \rangle . 
\end{equation} 
This net is finite~\cite{wass, xu}.  It is the net corresponding to the Wess-Zumino-Witten conformal field theory for the group manifold $G$, and to the vertex algebra associated to the affine Kac-Moody algebra $\widetilde{L \fg}$.
The nets associated to the loop groups of other Lie groups are also expected to be finite.  For example, combining the methods of Xu~\cite{xu} and the results of Toledano~\cite{toledano} should produce a proof that $\cN_{\mathrm{Spin}(2n),k}$ is finite.
\end{example}

In order to produce local field theories associated to these loop group conformal nets, we need to know that finite nets are dualizable.  In fact, finiteness provides a complete description of dualizable nets:
\newcounter{currenttheorem}
\setcounter{currenttheorem}{\value{theorem}}
\setcounter{theorem}{1}
\begin{theorem} %\label{thm-dualizable}
A conformal net $\cN \in \CN$ is dualizable if and only if it is finite, that is a finite direct sum of irreducible nets with finite $\mu$-index.
\end{theorem}
\begin{corollary} %\label{cor-lft}
For any finite conformal net $\cN$, there exists a framed local field theory $F_{\cN}\!:\! \FrBord_0^3 \!\ra\! \CN$ with $F_\cN(\,\pospnt) = \cN$.
\end{corollary}
\setcounter{theorem}{\value{currenttheorem}}
\nid Applied to the nets of Example~\ref{eg-lg}, this result provides a partial solution to the long-standing problem of finding local field theories associated to loop groups.

\subsection*{Finiteness of defects and sectors}
We can not only classify the set of framed local field theories with target conformal nets, as in Theorem~\ref{thm-dualizable}, but we can describe a 3-category of such local field theories.  A morphism between local field theories is, by definition, a local field theory on 3-manifolds equipped with partitioning codimension-1 domain walls.  A 2-morphism between such 1-morphisms is a local field theory on 3-manifolds with domain walls and partitioning codimension-1 submanifolds within the domain walls.  A 3-morphism is in turn a local field theory accommodating point partitions within the submanifolds within the domain walls.  Morphisms between net-valued local field theories are determined by dualizable defects, and more generally a 3-category of field theories is determined by explicitly identifying a sub-3-category $\CN^\fin \subset \CN$ that has duals~\cite{hopklurie}.  The morphisms and 2-morphisms of $\CN^\fin$ are characterized by the following finiteness properties.
\begin{definition}
An irreducible defect $\cD$ is \emph{finite} if the following conditions are satisfied: \begin{enumerate}
\item[1.] If $J$ and $J'$ are disjoint subintervals of $S^1$ containing respectively $i$ and $-i$, then the action of $\cD(J) \otimes^\alg \cD(J')$ on the sector $1_\cD$ extends to the spatial tensor product $\cD(J) \bar{\otimes} \cD(J')$. \vspace{2pt}
\item[2.] The bimodule
\begin{equation} \nn
_{\cD(J_1) \bar{\otimes} \cD(J_3)} (1_\cD)_{\cD(J_2)^{\op} \bar{\otimes} \cD(J_4)^{\op}}
\end{equation}
has finite statistical dimension.
\end{enumerate}
A defect is finite if it is a finite sum of finite irreducible defects.
\end{definition}
\begin{definition}
A sector $H$ between irreducible defects $\cD$ and $\cE$ is \emph{finite} if the bimodule $_{\cD(S^1_+)} H_{\cE(S^1_-)^{\op}}$ has finite statistical dimension.  In general a sector is finite if it is a finite sum of finite sectors between irreducible defects.
\end{definition}
\begin{theorem} \label{thm-cnfin}
Let $\CN^\fin$ denote the sub-symmetric monoidal 3-category of $\CN$ consisting of finite nets, finite defects, finite sectors, and all maps of finite sectors.  This 3-category $\CN^\fin$ has duals.
\end{theorem}
\nid In particular, $\CN^\fin$ encodes a 3-category of local field theories and their transformations.  It may be the case that $\CN^\fin$ is the maximal sub-3-category of $\CN$ that has duals, in which case $\CN^\fin$ is a full model for the collection of conformal net-valued local field theories, or, viewed alternately, for the collection of finite CFTs.

\subsection*{Proof of Theorem~\ref{thm-cnfin}}
Let $S^1_l$ and $S^1_r$ denote the halves of the standard circle with nonpositive and nonnegative real coordinate, respectively.  (1) The dual $\cN^\vee$ of a finite net $\cN$ is defined by
\begin{equation} \nn
\cN^\vee(I) = \cN(I)^\op,
\end{equation}
and is itself finite.  The counit defect $\cD_v: \cN \otimes \cN^\vee \ra 1$ of this duality is defined on a genuinely bicolored interval by
\begin{equation} \nn
\cD_v(I) = \cN(I_w \cup S^1_r \cup I_w).
\end{equation}
The unit defect $\cD_u: 1 \ra \cN^\vee \otimes \cN$ is similarly defined as
\begin{equation} \nn
\cD_u(I) = \cN(I_b \cup S^1_l \cup I_b).
\end{equation}
There is an invertible sector between the fusion $(1 \otimes \cD_u) \cdast (\cD_v \otimes 1)$ and the identity, as required.  (2) By the split property and the finiteness of the net $\cN$, the defects $\cD_v$ and $\cD_u$ are themselves finite.  (3) For a finite defect $\cD : \cN \ra \cM$ there is a finite defect $\cD^\dag : \cM \ra \cN$ that is both a left and a right adjoint to $\cD$.  This adjoint defect is defined by 
\begin{equation} \nn
\cD^\dag(I) = \cD(I^\rev),
\end{equation}
where $I^\rev$ is the color reversed interval, that is $I^\rev = I$ with $I^\rev_w = I_b$ and $I^\rev_b = I_w$.  The underlying Hilbert space of the counit $C: \cD \cdast \cD^\dag \ra 1_\cN$ is $L^2(\cD(S^1_+))$.  Because of the finiteness conditions on the defect $\cD$, there are equivariant actions of $(\cD \cdast \cD^\dag)(I)$ and of $\cN(J)$ on $L^2(\cD(S^1_+))$, for $-i \nin I \subset S^1$ and $i \nin J \subset S^1$; these actions give $L^2(\cD(S^1_+))$ the structure of a sector from $\cD \cdast \cD^\dag$ to $1_\cN$.  The unit of the adjunction $\cD \dashv \cD^\dag$ and the counit and unit of the adjunction $\cD^\dag \dashv \cD$ are similar.  (4)  By the second finiteness condition on the defect $\cD$, the units and counits of the adjunctions $\cD \dashv \cD^\dag$ and $\cD^\dag \dashv \cD$ are finite sectors.  (5) For a finite sector $H : \cD \ra \cE$ there is a finite sector $\overline{H} : \cE \ra \cD$ that is both a left and a right adjoint to $H$.  For $i \nin I$, the action of $d \in \cD(I)$ on $\overline{H}$ is defined by $d \overline{v} = \overline{d^* v}$, and for $-i \nin J$ the action of $e \in \cE(J)$ on $\overline{H}$ is defined by $e \overline{v} = \overline{e^* v}$.
% [End proof.]

\subsection*{Structure group, tensor categories, and examples}

The local field theories provided by Corollary~\ref{cor-lft} are a priori framed local field theories.  It turns out, however, that they always lift to stably framed theories:
\begin{theorem} \label{thm-sf}
Any framed local field theory $F: \FrBord_0^3 \ra \CN$ with target conformal nets lifts, along the natural map $\FrBord_0^3 \ra \SFBord_0^3$ from framed to stably framed bordism, to a stably framed local field theory $\widetilde{F}: \SFBord_0^3 \ra \CN$.
\end{theorem}
\nid This result is particularly relevant because in this range of dimensions, a stably framed structure is equivalent to a string structure, which in turn is a slight refinement of the classical $p_1$-structures (aka riggings) prevalent in the literature on 3-dimensional topological quantum field theory.  Thus net-valued local field theories make direct contact with the existing formulations of 3-dimensional field theory.

We mention a conjectural relationship of conformal nets to another double deloop of vector spaces, namely tensor categories.  Associated to a net $\cN$ there is a braided monoidal category
\begin{equation} \nn
\Rep(\cN) := \Hom_\CN(1_\cN, 1_\cN),
\end{equation}
and this association should extend to a trifunctor $\CN \xra{\Rep} \TC$.  Given such a trifunctor, Theorem~\ref{thm-cnfin} would imply that the braided monoidal category $\Rep(\cN)$ of representations of any finite net $\cN$ is a dualizable tensor category.

We conclude with a few examples.
\begin{example}[Loop group fusion categories] \label{eg-lgfc}
We have already mentioned the dualizable net $\cN_{G,k}$ associated to a central extension of the loop group of $G=SU(n)$, and its associated field theory
\begin{equation} \nn
F_{\cN_{G,k}} : \FrBord_0^3 \ra \CN.
\end{equation}
The representation category $\Rep(\cN_{G,k})$ of the net $\cN_{G,k}$ is a fusion category thought to be isomorphic to the category $\Rep(LG,k)$ of representations of the loop group $LG$ at level $k$.  
This fusion category $\Rep(LG,k)$ should be a dualizable tensor category, with its associated $\TC$-valued local field theory
\begin{equation} \nn
F_{\Rep(LG,k)} : \FrBord_0^3 \ra \TC.
\end{equation}
The net-valued theory contains more information than the tensor-category-valued theory---for example, for a closed framed 3-manifold $M$, the associated 3-manifold invariants $F_{\Rep(LG,k)}(M) \in \CC$ and $F_{\cN_{G,k}}(M) \in \CC$ should satisfy
\begin{equation} \nn
F_{\Rep(LG,k)}(M) = | F_{\cN_{G,k}}(M) |^2.
\end{equation}
\end{example}

\begin{example}[Moonshine net] \label{eg-mn}
Kawahigashi and Longo constructed a finite conformal net $\Moon$, the Moonshine net, as a twisted orbifold construction on the net associated to the Leech lattice~\cite{kawalongo}.  By Corollary~\ref{cor-lft} there is an associated net-valued local field theory.  Kawahigashi and Longo prove that the automorphism group of the Moonshine net is the Monster group: $\Aut(\Moon) = \Monst$.  We therefore have a map from the classifying space of the Monster group into the space of local field theories.  One among many open questions is whether the Moonshine net is higher Morita equivalent to the $24^{\mathrm{th}}$ power of the Dirac fermion.
\end{example}

\begin{example}[Link and surface invariants]
As a consequence of Theorem~\ref{thm-cnfin}, all of the aforementioned net-valued local field theories generalize to local field theories with insertions.  In particular, given a finite net $\cN$, the associated local field theory can be evaluated on 3-manifolds with embedded links labeled by finite 2-morphisms $V: 1_\cN \ra 1_\cN$ and on 3-manifolds with embedded surfaces labeled by finite endomorphisms $\cD: \cN \ra \cN$.  The resulting invariants of 3-manifolds with embedded surfaces are new even for the loop group nets.  Furthermore, the field theories for different nets can be combined: there are invariants for manifolds decorated by any number of interacting nets, defect-labeled walls, and sector-labeled embedded links.
\end{example}

\subsubsection*{Acknowledgments}
We thank Michael M\"uger, Michael Hopkins, Jacob Lurie, Ulrich Pennig, Peter Teichner, and Stephan Stolz for informative discussions.  The second author was supported by a Miller Research Fellowship.

\end{document}